\documentclass[conference]{IEEEtran}
\ifCLASSINFOpdf
  \usepackage[pdftex]{graphicx}
\else
\usepackage[dvips]{graphicx}
\usepackage{epsfig}
\fi

\usepackage{amsmath, amssymb, amsthm}	
\usepackage{color}

\newcommand{\R}{{I\!\!R}}
\newcommand{\N}{{I\!\!N}}

\begin{document}
%
\title{Splitting Approach for Solving Multi-Component Transport Models with Maxwell-Stefan-Diffusion}

\author{\IEEEauthorblockN{J\"urgen Geiser}
\IEEEauthorblockA{Ruhr University of Bochum, \\
The Institute of Theoretical Electrical Engineering, \\
Universit\"atsstrasse 150, D-44801 Bochum, Germany \\
Email: juergen.geiser@ruhr-uni-bochum.de}}


%


\maketitle

\begin{abstract}

In this paper, we present splitting algorithms to 
solve multicomponent transport models with Maxwell-Stefan-diffusion approaches.
The multicomponent models are related to transport problems, while we consider
plasma processes, in which the local thermodynamic equilibrium and weakly ionized plasma-mixture
models are given. Such processes are used for medical and technical applications.
These multi-component transport modelling equations are related to
convection-diffusion-reactions equations, which are wel-known in transport
processes. The multicomponent transport models can be derived from the microscopic
multi-component Boltzmann equations with averaging quantities and
leads into the macroscopic mass, momentum and energy equations, which are nearly Navier-Stokes-like equations.
An additional extension of the multicomponent diffusion term is based on the Maxwell-Stefan approach.
Such an approach allows to derive a diffusivity matrix, while the molecular force are balanced and relate all individual species velocities.
The Maxwell-Stefan diffusion approach is nonlinear, while we need additional
balancing equations. Such additional nonlinear equations are solved with iterative schemes.
We concentrate on solving the mass conservation equations of the
Navier-Stokes-like equations. Here, we consider splitting approaches with
non-iterative and iterative schemes. The splitting approaches are effective
methods to decompose delicate multicomponent transport models, while
the different operators in the transport models can be solved with
optimal numerical solvers.
We discuss the benefits of the decomposition into the convection, diffusion
and reaction parts, which allows to use fast numerical solvers for each part.
Additional, we concentrate on the nonlinear parts of the multicomponent diffusion,
which can be effectively solved with iterative splitting approaches
In the numerical experiments, we see the benefit of combining iterative
splitting methods with nonlinear solver methods, while these methods can relax the nonlinear terms.
In the outview, we discuss the future investigation of the next steps in
our multicomponent diffusion approaches.

\end{abstract}


%
\IEEEpeerreviewmaketitle

\section{Introduction}

In this paper, we concentrate on applying splitting approaches
to multicomponent transport models with Maxwell-Stefan-Diffusion, \cite{sene06}.
The modelling equations can be derived from the linearized Boltzmann
equations with approximated collision-terms to Navier-Stokes-like equations,
while we apply the Chapman-Enskog expansion and averaging techniques concerning
particle density, particle flux and particle kinetic energy, see \cite{giovang1999} and \cite{giovang2010}.

In our modelling problem, we consider a simplified plasma model, which considers only
heavy particles, which can be modelled as following: \\

The distribution function of the heavy particles are given as $f_i({\bf x}, {\bf c}_i, t)$, while ${\bf x}$ is the three-dimensional
spatial coordinate, ${\bf c}_i$ is the velocity of the molecule and $t$ is the time.
The heavy-particle species distribution are given as: 
\begin{eqnarray}
&& {\cal D}_i(f_i) = {\cal S}_i(f_i) + {\cal C}_i(f_i) , \; i \in I 
\end{eqnarray}
where ${\cal S}_i(f_i)$ is the scattering source term, ${\cal C}_i(f_i)$ is the reactive source term and  ${\cal D}_i(f_i)$ is the differential operator, see \cite{giovang1999}.
By applying Chapman-Enskog expansion, we apply the first-order perturbed distribution function to the linearized Boltzmann equations
and obtain with the averaging quantities of the particle density the macroscopic equation (mass-conservation). The mass conservation is given as:
\begin{eqnarray}
\label{navier_1}
&& \frac{\partial \rho_i}{\partial t} + \nabla_{\bold x} \cdot ( \rho_i \bold v) +  \nabla_{\bold x} (\rho_i {\cal V}_i) = m_i \omega_i^0 , \; i \in I , 
\end{eqnarray}
while $\rho_i$ is the mass density of $i$, ${\bold v}$ is the mean velocity, ${\cal V}_i$ is the species diffusion velocities and $\omega_i^0$ is the zero-th order production rate of species $i$, see \cite{giovang1999}.

We apply the Maxwell-Stefan approach and consider 3 species, then the transport model of the species, see the derivation in \cite{boudin2012}, which is given as:
\begin{eqnarray}
\label{ord_0}
&& \partial_t \xi_i - \nabla \cdot ({\bf v} \xi_i) + \nabla \cdot N_i = S_i , \; 1 \le i \le 3 , \\
&& \sum_{j=1}^3 N_j = 0 , \\
&& \frac{\xi_2 N_1 - \xi_1 N_2}{D_{12}} + \frac{\xi_3 N_1 - \xi_1 N_3}{D_{13}} =  -  \nabla \xi_1 , \\
 && \frac{\xi_1 N_2 - \xi_2 N_1}{D_{12}} + \frac{\xi_3 N_2 - \xi_2 N_3}{D_{23}} =  -  \nabla \xi_2 ,
\end{eqnarray}
where $\xi_i$ are the mole fractions and $N_i$ is the molar flux of
species $i$, see \cite{bothe2011} and \cite{boudin2012}.
The velocity ${\bf v}$ is given by the Navier-Stokes equation and we assume to deal with a 
first simplified model with a constant velocity.
Furthermore, the kinetic term or reaction term $S_i$ is given as:
\begin{eqnarray}
\label{ord_1_1}
&& S_i =  \sum_{j=1}^3 \lambda_{i,j} \xi_j  ,
\end{eqnarray}
where $\lambda_{i,j}$ are the reaction-rates.
The domain is given as $\Omega \in \R^d, d \in \N^+$ with $\xi_i \in C^2$.

\section{Methods}

We apply Operator splitting techniques to decompose the delicate full
differential equations, see \cite{geiser2008} and \cite{geiser2011}.
We decompose into a diffusion, a reaction and a convection part, see \cite{geiser2016} and \cite{geiser2016_1}. We apply the following 
splitting approach to our problem, we compute $n = 1, \ldots, N$, $t_0, t_1, \ldots, t_n$ time-steps:

\begin{itemize}

\item
The first step is given as (Diffusion step):
\begin{eqnarray}
\label{ord_0}
&& \hspace{-0.5cm} \partial_t \tilde{\xi}_i + \nabla \cdot N_i = 0 , \; 1 \le i \le 3 , \;\mbox{for} \; t \in [t^n, t^{n+1}], \\
&& \sum_{j=1}^3 N_j = 0 , \\
&& \frac{\tilde{\xi}_2 N_1 - \tilde{\xi}_1 N_2}{D_{12}} + \frac{\tilde{\xi}_3 N_1 - \xi_1 N_3}{D_{13}} =  -  \nabla \xi_1 , \\
 && \frac{\tilde{\xi}_1 N_2 - \tilde{\xi}_2 N_1}{D_{12}} + \frac{\tilde{\xi}_3 N_2 - \tilde{\xi}_2 N_3}{D_{23}} =  -  \nabla \tilde{\xi}_2 , \\
&& \tilde{\xi}_i(t^n) = \xi_i(t^n) , \; i = 1,2,3, \; t \in [t^n, t^{n+1}],
\end{eqnarray}
\item the next step is given as (Reaction step): 
\begin{eqnarray}
\label{ord_0}
&& \partial_t \hat{\xi}_i = S_i , \; 1 \le i \le 3 , \mbox{for} \; t \in [t^n, t^{n+1}], \\
&& \hat{\xi}_i(t^n) = \tilde{\xi}_i(t^{n+1}) , \; i = 1,2,3 .
\end{eqnarray}

\item and the next step is given as (Convection step): 
\begin{eqnarray}
\label{ord_0}
&& \hspace{-0.5cm} \partial_t \xi_i =  \nabla \cdot ({\bf v} \xi_i) , \; 1 \le i \le 3 , \mbox{for} \; t \in [t^n, t^{n+1}], \\
&& \xi_i(t^n) = \hat{\xi}_i(t^{n+1}) , \; i = 1,2,3 .
\end{eqnarray}

\end{itemize}

In the following section, we will discuss the results of the methods in an application.

\section{Results}

We apply model of a hydrogen plasma, see \cite{sene06}, which
we have simplified in a mass-transport model. We deal with the
species $H, H_2, H_2^+$, which we are heavy particles.
We take into account the following reactions, which are
given as:
\begin{eqnarray}
H_2 + e \; \underleftrightarrow{\lambda_1} \; H_2^+ + 2 e , \\
H_2 + e \; \underleftrightarrow{\lambda_2} \; 2 H + e ,
\end{eqnarray}
where the electron temperature is given as $T_e = 17400 \; [K]$ 
and the gas temperature values remain constant $T_h = 600 \; [K]$. \\

The simplified three component system is given as:
\begin{eqnarray}
\label{ord_0}
&& \partial_t \xi_i - v \partial_x \xi_i + \partial_x N_i = 0 , \; 1 \le i \le 3 , \\
&& \sum_{j=1}^3 N_j = 0 , \\
&& \frac{\xi_2 N_1 - \xi_1 N_2}{D_{12}} + \frac{\xi_3 N_1 - \xi_1 N_3}{D_{13}} =  -  \partial_x \xi_1 , \\
 && \frac{\xi_1 N_2 - \xi_2 N_1}{D_{12}} + \frac{\xi_3 N_2 - \xi_2 N_3}{D_{23}} =  -  \partial_x \xi_2 ,
\end{eqnarray}
where the domain is given as $\Omega \in \R^d, d \in \N^+$ with $\xi_i \in C^2$.

The parameters and the initial and boundary conditions used as following:
\begin{itemize}
\item $v = 0.01$ (means we have a small convection instead of the diffusion)
\item $D_{12} = D_{13} = 0.833$ (means $\alpha = 0$) and $D_{23} = 0.168$ (uphill diffusion, semi-degenerated Duncan and Toor experiment)
\item $D_{12} = 0.0833, D_{13} = 0.680$ and $D_{23} = 0.168$ (asymptotic behavior, Duncan and Toor experiment)
\item $J = 140$ (spatial grid points)
\item  The time-step-restriction for the explicit method is given as: \\
 $\Delta t \le (\Delta x)^2 \max \{ \frac{1}{2 \{D_{12}, D_{13}, D_{23}\}} \}$
\item The spatial domain is $\Omega = [0, 1]$, the time-domain $[0, T] = [0, 1]$
\item The initial conditions are:
\begin{eqnarray}
  \label{init}
  \xi_1^{in}(x) = \left\{ \begin{array}{l l}
0.8 & \mbox{if} \; 0 \le x < 0.25 , \\
1.6 (0.75 - x) & \mbox{if} \; 0.25 \le x < 0.75 , \\
0.0 & \mbox{if} \; 0.75 \le x \le 1.0 , 
\end{array} \right. ,
\end{eqnarray}
\begin{eqnarray}
&& \xi_2^{in}(x) = 0.2 , \; \mbox{for all} \; x \in \Omega = [0,1] ,
\end{eqnarray}
\item The boundary conditions are of no-flux type:
\begin{eqnarray}
\label{init}
&& N_1 = N_2 = N_3 = 0 , \mbox{on} \; \partial \Omega \times [0,1] ,
\end{eqnarray}
\end{itemize}

\begin{figure}[ht]
\begin{center}  
\includegraphics[width=5.0cm,angle=-0]{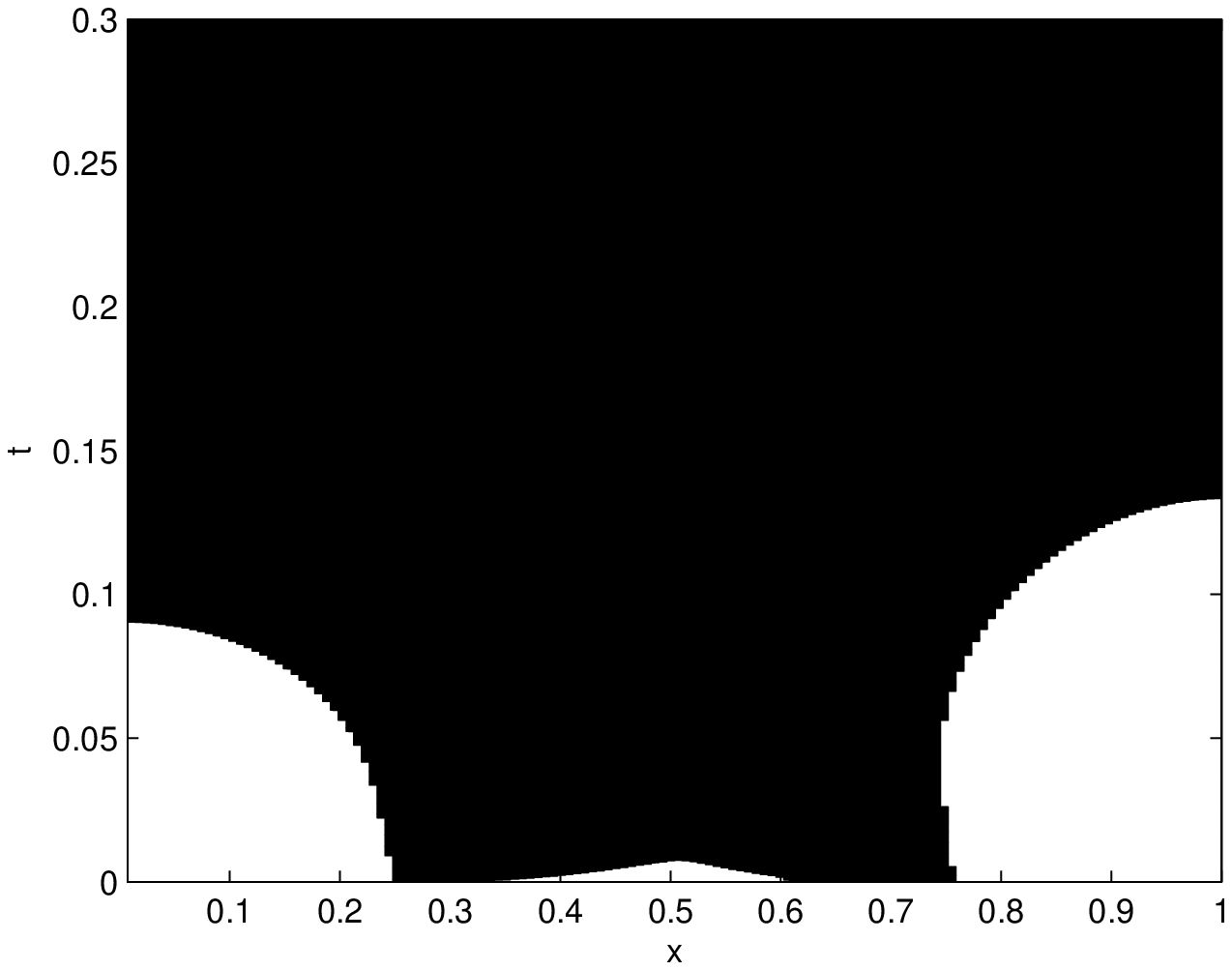} 
\includegraphics[width=5.0cm,angle=-0]{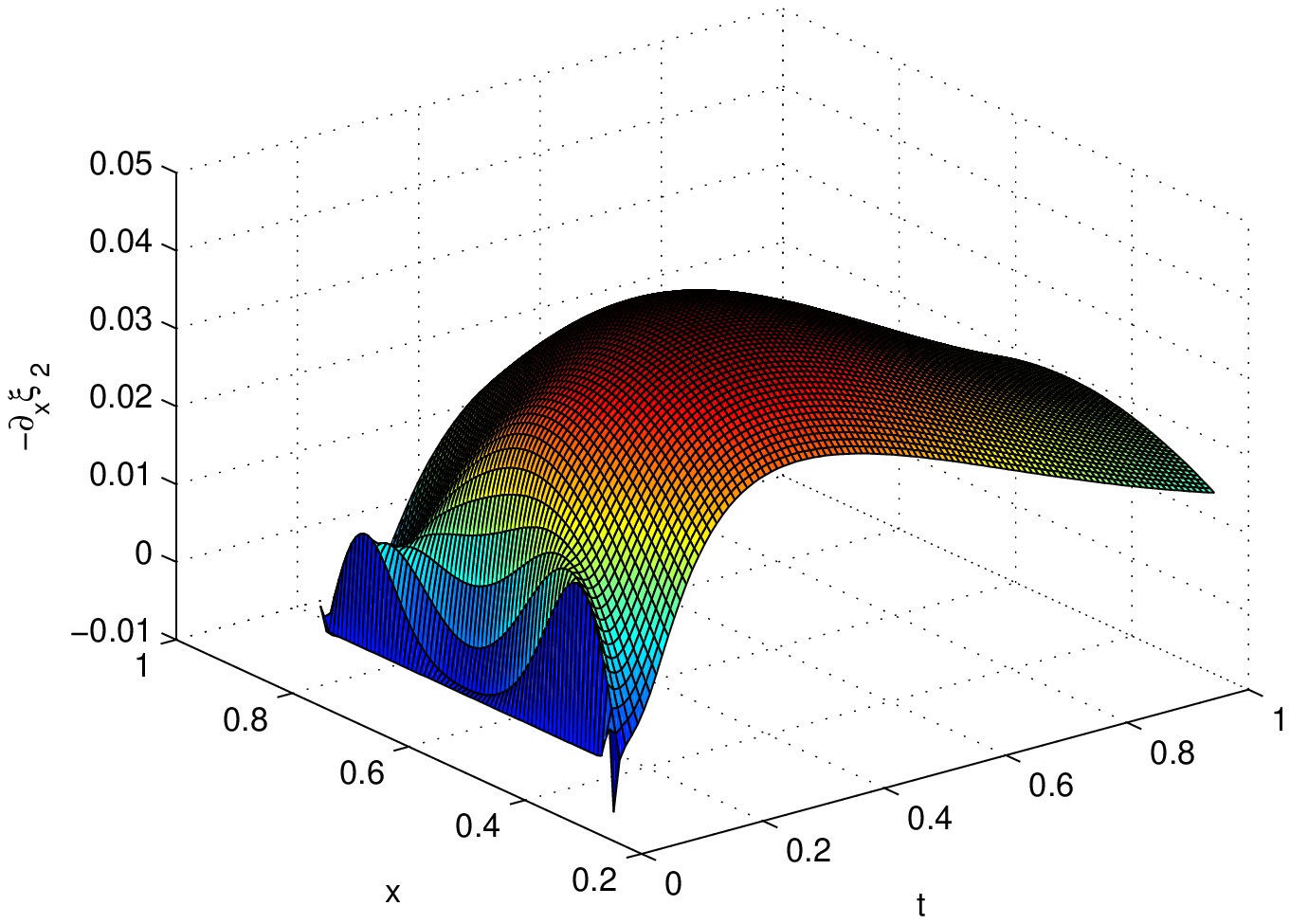} 
\end{center}
\caption{\label{hydrogen_3_2} The upper figure presents the result in the space time region, the lower figure presents the the 3D plot of the second component and.}
\end{figure}

The numerical solutions of the three hydrogen plasma in experiment 3 with the
uphill diffusion \ref{hydrogen_3_2}.

\section{Conclusions and Contributions}

We present the coupled model for a multi-component transport model simplified plasma processes. The Maxwell-Stefan diffusion
approach is considered and solved with additional iterative methods.
The nonlinear partial differential equations are splitted into a convection-, diffusion- and reaction part
and solved separately with optimal spatial discretization and time-integrator methods.
The numerical algorithms are presented and their numerical
convergences can be shown, see \cite{geiser2016_1}.
Although iterative splitting methods have the benefit of relaxing the nonlinearities and
can be used additional with their splitting approaches.
They are more accurate than noniterative splitting approaches.
The implicit behavior of iterative methods allows larger
time-steps to be used and they could accelerate the solver process.
In the future we aim to study the numerical analysis of the different combined schemes
and enlarge to real-life examples with more species and additional momentum- and energy equations.

\end{document}